\documentclass[aop,MSNbibl,dvips]{arximspdf}

\doi{10.1214/12-AOP744} 
\volume{41}
\issue{2}
\pubyear{2013}
\firstpage{671}
\lastpage{698}

\makeatletter
\newproclaim{example}{Example}[section]
\newtheorem{theorem}{Theorem}[section]
\newtheorem{lemma}{Lemma}[section]
\newtheorem{corollary}{Corollary}[section]
\newproclaim{remark}{Remark}[section]

\def\ov{\overline}

\def\al{\alpha}
\def\bb{\beta}
\def\Ga{\Gamma}
\def\De{\Delta}
\def\ep{\varepsilon}
\def\la{\lambda}
\def\ka{\kappa}
\def\Om{\Omega}
\def\si{\sigma}
\def\th{\theta}
\def\ze{\zeta}

\def\wt{\widetilde}
\def\rar{\rightarrow}
\def\stl{\stackrel{\LL}{=}}
\def\cd{\cdot}

\def\FF{\mathcal{F}}
\def\LL{\mathcal{L}}
\def\SS{\mathcal{S}}

\makeatother

\begin{document}
\begin{frontmatter}

\title{A sufficient condition for the continuity of permanental
processes with applications to local times of Markov processes\thanksref{T2}}
\runtitle{\hspace*{-10pt}Permanental processes and local times of Markov processes}

\begin{aug}
\author[A]{\fnms{Michael B.} \snm{Marcus}\corref{}\ead[label=e1]{mbmarcus@optonline.net}\ead[label=u1,url]{http://home.earthlink.net/\textasciitilde mbmarcus}}
\and
\author[B]{\fnms{Jay} \snm{Rosen}\ead[label=e2]{jrosen30@optimum.net}\ead[label=u2,url]{http://www.math.csi.cuny.edu/\textasciitilde rosen}}
\thankstext{T2}{Supported by grants from the NSF and PSC-CUNY.}
\runauthor{M. B. Marcus and J. Rosen}
\affiliation{City College of New York and College of Staten Island}
\address[A]{Department of Mathematics\\
City College of New York\\
NAC 6-291C\\
New York, New York 10031\\
USA\\
\printead{e1}\\
\printead{u1}}

\address[B]{Department of Mathematics\\
College of Staten Island\\
2800 Victory Blvd.\\
Staten Island, New York 10314\\
USA\\
\printead{e2}\\
\printead{u2}}

\end{aug}

\received{\smonth{5} \syear{2010}}
\revised{\smonth{11} \syear{2011}}

%
\begin{abstract}
We provide a sufficient condition for the continuity of real valued
permanental processes.
When applied to the subclass of permanental processes
which consists of squares of Gaussian processes, we obtain the
sufficient condition for continuity which is also known to be
necessary. Using an isomorphism theorem of Eisenbaum and Kaspi which
relates Markov local times and permanental processes, we obtain a
general sufficient condition for the joint continuity of local times.
\end{abstract}

%


\end{frontmatter}

\section{Introduction} \label{sec-1}

Let $T$ be an index set and $\{G(x),x\in T\}$ be a mean zero Gaussian
process with covariance $u(x,y)$, $x,y\in T$. It is remarkable that for
certain Gaussian processes, called associated processes, the process
$G^{2}=\{G^{2}(x),x\in T\}$ is closely related to the local times of a
strongly symmetric Borel right process with zero potential density
$u(x,y)$. This connection was first noted in the Dynkin Isomorphism
theorem~\cite{D83,D84} and has been studied by several probabilists,
including the authors and Eisenbaum and Kaspi. Our book~\cite{book}
presents several results about local times that are obtained using this
relationship.

The process $G^{2}$ can be defined by the Laplace transform of its
finite joint distributions
%
\begin{equation}
E \Biggl(\exp \Biggl(-\frac{1}{2}\sum_{i=1}^{n}
\al_{i}G^{2} (x_{i}) \Biggr) \Biggr)=
\frac{1}{|I+\al U |^{1/2}}\label{18q}
\end{equation}
for all $x_{1},\ldots,x_{n}$ in $T$, where $I$ is the $n\times n$
identity matrix, $\al$ is the diagonal matrix with $ (\al_{i,i}=\al_{i} )$, $\al_{i}\in R_{+}$ and $U=\{u(x_{i}, x_{j})\}
$ is an $n\times n$ matrix, that is symmetric and positive definite.

In 1997, Vere-Jones~\cite{VJ} introduced
the permanental process $\th:=\{\th_{x}, x\in T\} $, which is a real
valued positive stochastic process with finite joint distributions that satisfy
%
\begin{equation}
E \Biggl(\exp \Biggl(-\frac{1}{2}\sum_{i=1}^{n}
\al_{i}\th_{x_{i}} \Biggr) \Biggr)=\frac{1}{|I+\al\Gamma|^{\bb}}\label{18},
\end{equation}
where $\Ga=\{\Ga(x_{i}, x_{j})\}_{i,j=1}^{n}$ is an $n\times n$
matrix and $\bb>0$. (It would be better to refer to $\th$ as a $\bb
$-permanental process.) In this paper, in analogy with (\ref{18q}),
we consider these processes only for $\bb=1/2$ and refer to them as
permanental processes. The generalization here is that $\Ga$ need not
be symmetric or positive definite.

Even in (\ref{18q}), the matrix $U$ is not unique. The determinant
%
\begin{equation}
|I+\al U |=|I+\al MUM |
\end{equation}
for any signature matrix $M$. (A signature matrix is a diagonal matrix
with entries $\pm1$.)

The nonuniqueness is even more evident in (\ref{18}). If $D$ is any
diagonal matrix with nonzero entries, we have
%
\begin{equation}
|I+\al\Ga|=\bigl|I+\al D^{-1}\Ga D\bigr |=\bigl|I+\al D^{-1}
\Ga^{T} D \bigr|.
\end{equation}
For a very large class of irreducible matrices $\Ga$, it is known that
these are the only sources of nonuniqueness; see~\cite{L}.
On the other hand, in certain extreme cases, for example, if $\Ga_{1}$
and $\Ga_{2}$ are $n\times n$ matrices with the same diagonal elements
and all zeros below the diagonal, then $|I+\al\Gamma_{1} |=|I+\al
\Gamma_{2} |$. For this reason we refer to a matrix $\Ga$ for which
(\ref{18}) holds as a kernel of $\th$ (rather than as the kernel of
$\th$).

When $\Ga$ is not symmetric and positive definite, it is not at all
clear what kernels~$\Ga$ allow an expression of the form (\ref{18}).
(In~\cite{VJ} necessary and sufficient conditions on $\Ga$ for (\ref
{18}) to hold are given, but they are very difficult to verify. There
are very few concrete examples of permanental processes in~\cite{VJ}.)

It follows from the results in~\cite{VJ} that a sufficient condition
for (\ref{18}) to hold is that all the real nonzero eigenvalues of
$\Ga$ are positive and that $r\Ga(I+r\Ga)^{-1}$ has only
nonnegative entries for all $r>0$. \label{page-3} In~\cite{EK},
Eisenbaum and Kaspi note that this is the case when $\Ga(x,y)$,
$x,y\in T$, is the potential density of a transient Markov process on
$T$. This enables them to find a Dynkin-type isomorphism for the local
times of Markov processes that are not necessarily symmetric, in which
the role of $G^{2}$ is taken by the permanental process $\th$.

Both Eisenbaum and Kaspi have asked us if we could find necessary and
sufficient conditions for the continuity and boundedness of permanental
processes. In this paper we give a sufficient condition for the
continuity of permanental processes.
When applied to the subclass of permanental processes
which consists of squares of Gaussian processes, it is, effectively,
the sufficient condition for continuity which is also known to be
necessary. We use our sufficient condition for the continuity of
permanental processes and an isomorphism theorem for permanental
processes given by Eisenbaum and Kaspi in~\cite{EK}, Theorem 3.2, to
extend a sufficient condition they obtain in~\cite{EK1}, Theorem 1.1,
for the continuity of local times of Markov processes, to a larger
class of Markov processes.

In Section \ref{sec-3} we review several properties of
permanental processes. In particular, a key property of permanental
processes is that $\Ga(x,x)\geq0$ and
%
\begin{equation}
0\le\Ga(x , y)\Ga(y, x )\le\Ga(x,x) \Ga(y,y)\qquad \forall x,y\in T. \label{15aa}
\end{equation}
This allows us to define
%
\begin{equation}
d(x,y) =4\sqrt{ 2/3} \bigl(\Ga(x,x)+\Ga(y,y)-2 \bigl(\Ga(x , y)\Ga (y, x )
\bigr)^{1/2} \bigr)^{1/2}.\label{12}
\end{equation}
Let $D=\sup_{s,t\in T} \,d (s,t)$. $D$ is called the $ d $ diameter of $T$.

Let $(T,\rho)$ be a metric or pseudometric space. Let $B_{\rho}(t,u)$
denote the closed ball in $(T,\rho)$ with radius $u$ and center $t$.
For any probability measure $\mu$ on $(T,\rho)$, we define
%
\begin{equation}
J_{T, \rho,\mu}( a) =\sup_{t\in T}\int_0^a
\biggl(\log\frac1{\mu \bigl(B_{\rho}(t,u)\bigr)} \biggr)^{1/2}
\,du.\label{tau}
\end{equation}
We occasionally omit some of the subscripts $T, \rho$ or $\mu$, if
they are clear from the context.

Whether or not $d (x,y)$ is a metric, or pseudometric on $T$, we can
define the sets $B_{ d}(s,u)=\{t\in T |  \,d (s,t)\leq u\}$.
We can then define $ J_{T, d,\mu}( a)$ as in (\ref{tau}), for any
probability measure $\mu$ for which the sets $B_{ d}(s,u)$ are measurable.

\begin{theorem}\label{theo-11a} Let $T$ be a separable topological
space, and
let $\mathcal{B}(T)$ denote it's Borel $\si$-algebra. Let $\th= \{
\th_x\dvtx  x\in T\}$ be a permanental process with kernel $\Ga$ with the
property that $\sup_{x\in T} \Ga(x,x)<\infty$. Assume that $ d
(x,y)$ is continuous on $T\times T$ and that there exists a probability
measure $\mu$ on $\mathcal{B}(T)$ such that
%
\begin{equation}
\lim_{\delta\to0}J_{ d}(\delta)=0.\label{18jv}
\end{equation}
Then there exists a version $\th'= \{\th'_x\dvtx  x\in T\}$ of $\th$ that
is bounded and continuous almost surely and satisfies
%
\begin{equation}
\lim_{\delta\to0} \mathop{\sup_{s,t\in T}}_{ d (s,t)\le\delta
}\frac{|\th'_{ s} -\th'_{t}|}{J_{ d} ( d(s,t) /2) }\le60 \Bigl(
\sup_{x\in T}\th'_x \Bigr)^{1/2}\qquad
\mbox{a.s.},\label{21sv}
\end{equation}
where in (\ref{21sv}) and in similar situations elsewhere in this
paper, we make the convention that $0/0=0$.
\end{theorem}

We show in Lemma~\ref{lem-33q} that when $\th$ is continuous on $T$
almost surely, then $ d (x,y)$ is continuous on $T\times T$. Therefore,
the condition in Theorem~\ref{theo-11a}, that $ d (x,y)$ is
continuous on $T\times T$, is perfectly reasonable. In particular, it
is implied by the continuity of $\Ga(x,y)$.

We say that a metric or pseudometric $d_{1}$ dominates $ d$ on $T$
if
%
\begin{equation}
d(x,y)\leq d_{1}(x,y) \qquad\forall x,y\in T.\label{dom}
\end{equation}
In Section~\ref{sec-4}, we give several natural metrics that dominate
$ d$.

\begin{corollary}\label{cor-11b} Let $\th= \{\th_x\dvtx  x\in T\}$ be a
permanental process with kernel~$\Ga$ satisfying $\sup_{x}\Ga
(x,x)<\infty$. Let $d$ be given by (\ref{12}), and let $d_{1}(x,y)$
be a metric or pseudo-metric on $T$ that dominates $ d (x,y)$ and is
such that $(T, d_{1})$ is separable and has finite diameter $D$.
Consider $T$ with the $d_{1}$ topology, that is, $(T,d_{1})$. Then
Theorem~\ref{theo-11a} holds with $d$ replaced by $d_{1}$.
\end{corollary}

In Section~\ref{sec-loc} we give a version of (\ref{21sv}) for $|\th'_{s}-\th'_{t_{0}}|$ for fixed $t_{0}\in T$, which provides a local
modulus of continuity for permanental processes.

Let $X =
(\Om, X_t,P^x
)$ be a transient Borel right process with state space $S$ and
$0$-potential density
$u(x,y)$. We assume that $S$ is a locally compact topological space,
and that $u(x,y)$ is continuous. These conditions imply that $X$ has
local times; see, for example,~\cite{book}, Theorem 3.6.3. It is shown
in~\cite{EK}, Theorem 3.1, that there exists a permanental process
$\th=\{\th_{y} ; y\in S\}$, with kernel $u(x,y)$, which they refer
to as the
permanental process associated with $X$.

In~\cite{EK}, Theorem 3.2,
an isomorphism theorem is given that relates the local times of $X$ and
$\th$.
In the next theorem, we use this isomorphism together with Theorem~\ref
{theo-11a} in this paper, to obtain a sufficient condition for the
joint continuity of the local times of $X$.
When applied to strongly symmetric Markov processes, we obtain the
sufficient condition for joint continuity, that is known to be
necessary; see~\cite{book}, Theorem 9.4.11. Applied to L\'{e}vy
processes, which need not be symmetric, we also obtain the sufficient
condition for the joint continuity of local times that is known to be
necessary; see~\cite{Bar84}.

As usual, we use $\ze$ to denote the death time of $X$.

\begin{theorem}\label{theo-jc} Let $S$ be a locally compact
topological space
with a countable base.
Let $X =
(\Om, X_t, P^x
)$ be a recurrent Borel right process with state space $S$ and
continuous, strictly positive
$1$-potential densities
$u^{1}(x,y)$. Define
$d(x,y)$ as in (\ref{12}) for the kernel $u^{1}(x,y)$.
Suppose that for every compact set $K\subseteq S$, we can find a
probability measure $\mu_{K}$ on $K$, such that
%
\begin{equation}
\lim_{\delta\to0} J_{K, d,\mu_{K}}( \delta) =0.\label{18jj}
\end{equation}
Then X has a jointly continuous local time
$\{L^y_{t} ; (y,t)\in S\times R_+\}$.\vadjust{\goodbreak}

Let $X$ be a transient Borel right process with state space $S$ and
continuous, strictly positive
$0$-potential densities
$u(x,y)$. If (\ref{18jj}) holds for every compact set $K\subseteq S$,
with $d(x,y)$ defined as in (\ref{12}) for the kernel $u (x,y)$, X
has a local time
$\{L^y_{t} ; (y,t)\in S\times R_+\}$ which is jointly continuous on
$S\times[0,\zeta)$.
\end{theorem}


In Theorem~\ref{theo-jc} we only get continuity of the local times of
transient processes on $S\times[0,\zeta)$. However, as it is pointed
out in~\cite{EK1}, if $X$ is transient, using an argument of Le Jan,
we can always find a recurrent process $Y$ such that $X$ is $Y$ killed
the first time it hits the cemetery state $\De$. Problematically, this
changes the potentials (see~\cite{DM}, (78.5)) and hence the condition
(\ref{18jj}). We leave it to the interested reader to work out the details.

It is interesting to place Theorem~\ref{theo-jc} in the history of
results on the joint continuity of local times of Markov processes. A
good discussion is given in~\cite{EK1}. We make a few comments here.
In~\cite{Bar84} Barlow gives necessary and sufficient condition for
the joint continuity of local times of L\'evy processes. Local times
are difficult to work with. He works hard to obtain many of their
properties. In~\cite{sip} we use the Dynkin Isomorphism theorem (DIT)
to obtain necessary and sufficient condition for the joint continuity
of local times of strongly symmetric Borel right processes, which,
obviously, includes symmetric L\'evy processes. Using the DIT enables
us to infer properties of local times from those of Gaussian processes.
These processes are well understood and easier to work with than local
times. Although the results in~\cite{sip} only give the results in
\cite{Bar84} for symmetric L\'evy processes, they apply to a much
larger class of symmetric Markov processes.

In~\cite{EK1}, Eisenbaum and Kaspi extend Barlow's approach to obtain
sufficient conditions for the joint continuity of local times of a
large class of recurrent Borel right processes and also give a modulus
of continuity for the local times. In Theorem~\ref{theo-jc}, using a
proof similar to the one in~\cite{sip}, we use Eisenbaum and Kaspi's
isomorphism theorem for permanental processes~\cite{EK}, Theorem 3.2,
to extend their results in~\cite{EK1}.
(In~\cite{EK1}, they require the existence of a Borel right dual
process. This is not needed in Theorem~\ref{theo-jc}. In Section \ref
{sec-comp} we show how to obtain~\cite{EK}, Theorem~3.2, from Theorem
\ref{theo-jc}.) We also obtain uniform and local moduli of continuity
for the local times.

\begin{theorem}\label{lt-um}
Under the assumptions of Theorem~\ref{theo-jc},
%
\begin{eqnarray}
&&\lim_{\delta\to0} \mathop{\sup_{x,y\in K}}_{d(x,y)\le\delta} \frac{|L^x_t-L^y_t|}{J_{K, d_{1},\mu_{K}}(d(x,y)/2)}\label{mod-1}
\nonumber
\\[-8pt]
\\[-8pt]
\nonumber
&&\qquad \le30\sup_{y\in K}\bigl(L^{y}_t
\bigr)^{1/2}\qquad \mbox{for almost all $t\in[0,\zeta)$ a.s.}
\end{eqnarray}
\end{theorem}

The local modulus of continuity for local times is given in Theorem
\ref{theo-62}.\vadjust{\goodbreak}

We thank Michel Talagrand for suggestions resulting in a
significant simplification of the proof of Theorem~\ref{theo-11a}.

\section{Some basic continuity theorems}\label{sec-2}

For $p\geq1$, let $
\psi_p(x)=
\exp(x^{p}) -1
$ and $L^{\psi_p}(\Om,\FF,P)$ denote the set of random
variables
$\xi\dvtx \Om\to R^{1}$ such that $E\psi_p ( |\xi|/c
)<\infty$ for some $c>0$.
$L^{\psi_p}(\Om,\FF,P)$ is a Banach space with norm given by
%
\begin{equation}
\|\xi\|_{\psi_p}=\inf \bigl\{c>0\dvtx E\psi_p \bigl( |\xi|/c \bigr)\le
1 \bigr\}.\label{14}
\end{equation}
We shall only be concerned with the cases $p=1$ and 2.

We obtain Theorem~\ref{theo-11a} with the help of the following basic
continuity theorems. They are, essentially, best possible sufficient
conditions for continuity and boundedness of Gaussian process. However,
it is well known that they hold for any stochastic process satisfying
certain conditions with respect to the Banach space~$L^{\psi_2}$.

\begin{theorem}\label{maj} Let $X= \{X(t)\dvtx  t\in T\}$ be a
stochastic process such that $X(t,\omega)\dvtx T\times\Om\mapsto
[-\infty,\infty] $ is $\mathcal{A}\times\mathcal{F}$ measurable
for some $\si$-algebra $\mathcal{A}$ on $T$. Suppose $X(t)
\in L^{\psi_2}(\Om,\FF,P)$, and let
%
\begin{equation}
\hat d(t,s):= \bigl\|X(t)-X(s)\bigr\|_{\psi_2}.\label{22ar}
\end{equation}
[Note that the balls $B_{\hat d}(s,u)$ are $\mathcal{A}$ measurable.]

Suppose that $(T,\hat d)$ has finite diameter $D$, and that there exists
a probability measure $\mu$ on $(T,\mathcal{A})$ such that
%
\begin{equation}
J_{\hat d}(D)<\infty. \label{210}
\end{equation}
Then there exists a version $X'=\{X'(t),t\in T\}$ of $X$ such that
%
\begin{equation}
E\sup_{t\in T}X'(t)\le C J_{\hat d}(D) \label{22}
\end{equation}
for some $C<\infty$. Furthermore for all $0<\delta\le D$,
%
\begin{equation}
\mathop{\sup_{s,t\in T}}_{\hat d(s,t)\le\delta} \bigl|X' ( s,\omega )-X'( t,\omega)\bigr|
\le2Z(\omega) J_{\hat d}(\delta)\label{212},
\end{equation}
almost surely, where
%
\begin{equation}
Z(\omega):=\inf \biggl\{\al>0\dvtx \int_{T}
\psi_{2}\bigl(\al^{-1}\bigl|X(t,\omega)\bigr|\bigr) \mu(dt) \leq1
\biggr\}
\end{equation}
and $\|Z\|_{\psi_2} \le
K $, where $K$ is a constant.

In particular, if
%
\begin{equation}
\lim_{\delta\to0}J_{\hat d}(\delta)=0,\label{25}
\end{equation}
$X'$ is uniformly continuous on $(T,\hat d)$ almost surely.
\end{theorem}

\begin{remark}\label{rem-orlicz}
Theorem~\ref{maj} is well known. It contains ideas that originated in
an important early
paper by Garcia, Rodemich and Rumsey Jr.~\cite{GRR}, and were
developed further by Preston~\cite{Preston1,Preston2} and Fernique
\cite{Fernique}. We present a generalization of it in~\cite{MRas}, Theorem
3.1. Unfortunately, the statement of~\cite{MRas}, Theorem~3.1,
makes it appear that (\ref{25}), in this paper, is required for (\ref
{212}), in this paper, to hold. This is not the case as one can see
from going through the proof of~\cite{MRas}, Theorem 3.1. However, an
easier way to see that (\ref{212}), in this paper, holds is to note
that it follows immediately from~\cite{book}, Theorem 6.3.3. Again,
unfortunately, the hypothesis of~\cite{book}, Theorem 6.3.3, requires
that $X$ is a Gaussian process. A~reading of the proof shows that it
actually only requires that $X(t)
\in L^{\psi_2}(\Om,P)$ and $\|X(t)-X(s)\|_{\psi_2} \le d(s,t) $ for
all $s,t \in
T$ where $d(s,t)$ is some metric; see also~\cite{KR}.
\end{remark}

The inequality in (\ref{212}) is not quite enough to give a best
possible uniform modulus of continuity for $X'$. Instead we use the
following lemma due to Heinkel~\cite{H}, Proposition 1.

\begin{lemma}\label{lem-22} Let $(T,\hat d)$ be a metric or
pseudo-metric space with finite diameter $D$ and $\mu$ be a
probability measure on $T$ with the property that $\mu(B_{\hat
d}(t,u))>0$ for all $t\in T$ and $u>0$.
Assume that (\ref{25}) holds.
Let $\{f(t),t\in T\}$ be continuous on $(T,\hat d)$, and set
%
\begin{equation}
\wt f(s,t)=\frac{f(s)-f(t)}{\hat d(s,t)}I_{\{(u,v)\dvtx \hat d(u,v)\ne0\}
}(s,t).\label{27}
\end{equation}
Then if
%
\begin{equation}
c_{\mu,T}(\wt f)=\int_ {T\times T}\psi_{2}\bigl(
\wt f(s,t)\bigr) \,d\mu(s) \,d\mu(t)<\infty,\label{c}
\end{equation}
we have that for all $x,y\in T$,
%
\begin{equation}
\bigl|f(x)-f(y)\bigr|\label{217}\le20 \sup_{t\in T}\int_{0}^{\hat
d(x,y)/2}
\biggl(\log \biggl(\frac{ c_{\mu,T}(\wt f)+1}{\mu^{2}(B_{\hat d}(t,u))} \biggr) \biggr)^{1/2} \,du.
\end{equation}
\end{lemma}

%


\begin{theorem}\label{theo-22} Under the hypotheses of Theorem \ref
{maj}, assume that (\ref{25}) holds. Then there exists a version
$X'=\{X'(t),t\in T\}$ of $X$ such that
%
\begin{equation}
\lim_{\delta\to0} \mathop{\sup_{s,t\in T}}_{\hat
d(s,t)\le
\delta}\frac{ |X'(s)-X'(t)|}{J_{\hat d}(\hat d(s,t)/2)}\le30 \qquad\mbox{a.s.}
\label{21sa}
\end{equation}
\end{theorem}

\begin{pf}
Assume first that we can find points $t_{1},\ldots, t_{n}$ such that
$\hat d(t_{i},t_{j})>0$ for all $i\neq j$. We can cover these points
with $n$ disjoint balls. Therefore the $\mu$ measure of one of these
balls must be less than or equal to $1/n$. Consequently, for all
$\delta>0$, sufficiently small
%
\begin{equation}
J_{T, \hat d,\mu}( \delta) \ge\delta(\log n)^{1/2}. \label{n215}\vadjust{\goodbreak}
\end{equation}
If $n$ is the maximal number of such points, then the $\sup$ on the
left-hand side of (\ref{21sa}) is zero for all $\delta$
sufficiently small. Here we use the fact that any other point $t\in T$
must satisfy $\hat d(t,t_{j})=0$ for some $ j$, and hence
$X_{t}=X_{t_{j}}$ a.s. by the definition of~$\hat d$, so that $\hat
d(t,t_{i})=\hat d(t_{i},t_{j})$ for all $i$. Thus (\ref{21sa}) is
trivially true.

If there is an infinite number of such points, it follows from (\ref
{n215}) that
%
\begin{equation}
\lim_{\delta\to0} \frac{ J_{\hat d}(\delta)}{\delta}=\infty \label{218a}.
\end{equation}

By Theorem~\ref{maj} we can assume that $X=\{X(t),t\in T\}$ is
continuous on $(T,\hat d)$ almost surely. Define $\wt X$ as in (\ref
{27}). Note that by Fubini's theorem
%
\begin{eqnarray}\quad
&& E \biggl(\int_ {T\times T}\psi_{2}\bigl( \wt X(s,t)
\bigr) \,d\mu(s) \,d\mu (t) \biggr)\label{213}
\nonumber
\\[-8pt]
\\[-8pt]
\nonumber
&&\qquad = E \biggl(\int_ {T\times T}\psi_{2} \biggl(
\frac{
X(t)-X(s)}{ \| X(t)-X(s)\|_{\psi_2}} \biggr)1_{\{0<\hat d(s,t)\}} \,d\mu(s) \,d\mu(t) \biggr)\leq1.
\end{eqnarray}
Consequently,
%
\begin{equation}
\int_ {T\times T}\psi_{2}\bigl( \wt X(s,t)\bigr) \,d
\mu(s) \,d\mu(t)<\infty\qquad \mbox{a.s.}
\end{equation}
Let $\Om'$ be the set of measure $1$ in the probability space for
which this is finite and for which $X(t, \omega)$ is continuous. For
each $\omega\in\Om'$,
%
\begin{equation}
c_{\mu,T}(\wt X)= \int_ {T\times T}\psi_{2}\bigl(
\wt X(s,t,\omega)\bigr) \,d\mu(s) \,d\mu(t) <\infty.
\end{equation}
To obtain (\ref{21sa}), we use (\ref{217}) with $f(\cd)$ replaced
by $\wt X(\cd)$. Note that the right-hand side of
(\ref{217})
%
\begin{equation}
\le10 \hat d(x,y) \bigl(\log \bigl( c_{\mu,T}(\wt f)+1 \bigr)
\bigr)^{1/2}+30 J_{\hat d}\bigl(\hat d(x,y)/2\bigr).\label{211}
\end{equation}
Using (\ref{218a}) allows us to simplify the denominator in (\ref
{21sa}).
\end{pf}

We get a result similar to (\ref{218a}) for the local modulus of
continuity, but it is more delicate. We take this up in Section~\ref{sec-loc}.

\section{\texorpdfstring{Proof of Theorem \protect\ref{theo-11a}}{Proof of Theorem 1.1}}\label{sec-3}

We begin with some observations about permanental processes.
It is noted in~\cite{VJ}, and immediately obvious from (\ref{18}),
that the univariate marginals
of a permanental process are squares of normal random variables. A key
observation used in the proof of Theorem~\ref{theo-11a}, which also
follows from (\ref{18}), is that the bivariate marginals of a
permanental process are squares of bivariate normal random variables.
We proceed to explain this.

For $n=2$, (\ref{18}) takes the form
%
\begin{eqnarray}\label{18a}
&& E \biggl(\exp \biggl(-\frac{1}{2} (\al_{1}
\th_{x }+\al_{2} \th_{y} ) \biggr) \biggr)
\nonumber
\\[-2pt]
&&\qquad =\frac{1}{|I+\al\Gamma|^{1/2}}
= \bigl(1+\al_{1}\Ga(x , x )+
\al_{2}\Ga(y, y)
\\[-2pt]
& &\hspace*{72pt}\qquad\quad{}
+\al_{1}\al_{2} \bigl(\Ga(x , x )\Ga
(y, y)-\Ga(x , y)\Ga(y, x ) \bigr) \bigr)^{-1/2}.\nonumber
\end{eqnarray}
Taking $\al_{1}=\al_{2}$ sufficiently large, this implies that
%
\begin{equation}
\Ga(x, x)\Ga(y, y)-\Ga(x, y)\Ga(y, x)\geq0. \label{32}
\end{equation}
If we set $\al_{2}=0$ in (\ref{18a}), we see that for any $x\in T$,
%
\begin{equation}
\Ga(x , x )\geq0\label{vj1}.
\end{equation}
In addition, by~\cite{VJ}, page 135, last line, for any pair $x ,y\in T$,
%
\begin{equation}
\Ga(x , y ) \Ga(y , x )\geq0.\label{vj2}
\end{equation}
It follows from (\ref{32})--(\ref{vj2}) that for any pair $x , y\in
T,$ the matrix
\[
\left[\matrix{ \Ga(x, x)& \bigl(\Ga(x, y)\Ga(y, x) \bigr)^{1/2}\vspace*{2pt}
\cr
\bigl(\Ga(x, y)\Ga(y, x) \bigr)^{1/2}&\Ga(y, y)}\right]
\]
is positive definite, so that we can construct a mean zero Gaussian
vector $\{G (x),G (y)\}$ with covariance matrix
%
\begin{equation}
E \bigl(G(x) G(y) \bigr) = \bigl(\Ga(x, y)\Ga(y, x) \bigr)^{1/2}.
\label{111a}
\end{equation}
Note that
%
\begin{equation}
\bigl(E\bigl(G(x)-G(y)\bigr)^{2} \bigr)^{1/2}=
\frac{\sqrt{3/2}}{4}\,d(x,y), \label{dd}
\end{equation}
defined in (\ref{12}).\vspace*{-2pt}

\begin{lemma}\label{lem-1} Suppose that $\th:=\{\th_{x}, x\in T\} $
is a permanental process
for $\Ga$ as given in (\ref{18}). Then for any pair $x , y$,
%
\begin{equation}
\{\th_{x},\th_{y}\}\stl\bigl\{G^{2}(x),G^{2}(y)
\bigr\},\label{fund}
\end{equation}
where $\{G (x),G (y)\}$ is a mean zero Gaussian random variable with
covariance matrix given by (\ref{111a}).\vspace*{-2pt}
\end{lemma}

\begin{pf} By (\ref{18a}) the Laplace transform of $\{\th_{x },\th_{y}\}
$ is the same as
the Laplace transform of $\{G^{2}(x ),G^{2}(y)\}$.\vspace*{-2pt}
\end{pf}

\begin{pf*}{Proof of Theorem~\ref{theo-11a}}
It follows from Lemma~\ref{lem-1} that
%
\begin{eqnarray}\label{ffun}
\hat d(x,y)&:=&\bigl\|\theta^{1/2}_{x}-\theta^{1/2}_{y}
\bigr\|_{\psi_2}=\bigl\| |G_{x}|-|G_{y}| \bigr\|_{\psi_2}
\nonumber
\\[-9pt]
\\[-9pt]
\nonumber
&\leq& \| G_{x}-G_{y} \|_{\psi_2}= d(x,y).\vadjust{\goodbreak}
\end{eqnarray}
Since
$d(x,y)$ is continuous, the metric $ \hat d(x,y)$ is also continuous.
Therefore, the separability of $T$ implies that $(T,\hat d)$ is a
separable metric space.
By~\cite{C}, Theorem~2, we may assume that
$\th^{1/2}=\{\th^{1/2}_x,x\in T\}$ is measurable with respect to
$(T,\hat d)$. (More explicitly, measurability means that $\th^{1/2}_{x}(\omega)\dvtx T\times\Om\mapsto[0,\infty] $ is $ \mathcal
{B}(T,\hat d)\times\mathcal{F}$ measurable.)

By
Theorem~\ref{maj} with $X=\th^{1/2}$ and $ \mathcal{A} = \mathcal
{B}(T,\hat d)$, we see that if there exists a probability measure $\mu
$ on $(T,\hat d)$ such that
%
\begin{equation}
\lim_{\delta\to0}J_{\hat d}(\delta)=0,\label{ffun25}
\end{equation}
then there exists a version $X'=\{X'(t),t\in T\}$ of $X$ such that
$X'$ is bounded and uniformly continuous on $(T,\hat d)$ almost surely.

By assumption, there exists a probability measure $\mu$ on $\mathcal
{B}(T)$ such that
%
\begin{equation}
\lim_{\delta\to0}J_{ d}(\delta)=0.\label{ffun18jv}
\end{equation}
Since $ \hat d(x,y)$ is continuous, $\mathcal{B}(T,\hat d)\subseteq
\mathcal{B}(T)$. Hence we can restrict $\mu$ to be a probability
measure on $(T,\hat d)$, and it follows from (\ref{ffun18jv}) and
(\ref{ffun}) that (\ref{ffun25}) holds. Thus we obtain a version
$X'$ which is bounded and continuous on $(T,\hat d)$, and using again
the continuity of $ \hat d(x,y)$, this implies continuity on~$T$.

Similarly, it follows from Theorem~\ref{theo-22} with $X=\th^{1/2}$ that
%
\begin{equation}
\lim_{\delta\to0} \mathop{\sup_{x,y\in T}}_{ d(x,y)\le\delta
}\frac{ |\theta^{1/2}_{x}-\theta^{1/2}_{y}|}{J_{ d}(d(x,y)/2)}\le 30.\label{21swb}
\end{equation}
Using the inequality
%
\begin{equation}
|\theta_{x}-\theta_{y}|\leq\bigl|\theta^{1/2}_{x}-
\theta^{1/2}_{y}\bigr| 2\sup_{z}\theta^{1/2}_{z},
\end{equation}
we get (\ref{21sv}).
\end{pf*}


\begin{pf*}{Proof of Corollary~\ref{cor-11b}}
By (\ref{ffun}),
$\hat d(x,y)\leq d_{1}(x,y)$. Consequently, the proof of Corollary
\ref{cor-11b} follows immediately from the proof of Theorem~\ref{theo-11a}.
\end{pf*}

\begin{lemma} \label{lem-33q} When $\th$ is continuous on $T$ almost
surely, $ d (x,y)$ is continuous on $T\times T$.
\end{lemma}

\begin{pf} By Lemma~\ref{lem-1},
%
\begin{equation}
E(\th_{ x})=\Ga(x,x) \quad\mbox{and}\quad \operatorname{cov}\{\th_{
x},
\th_{y}\}=2\Ga(x,y)\Ga(y,x).\label{42v}
\end{equation}
In addition, since the univariate marginals of $\th$ are the squares
of Gaussian random variables, $\th_{x}$ and $\th_{y}$ are locally
uniformly bounded in any $L^{p}$ space.
\end{pf}

\begin{remark} Theorem~\ref{maj} can be used to obtain more
information about $\th$. For example, a very minor modification of the
proof of Theorem~\ref{theo-11a} shows that when
%
\begin{equation}
J_{ d}(D)<\infty, \label{210a}
\end{equation}
there exists a version $X'=\{X'(t),t\in T\}$ of $X$ such that
%
\begin{equation}
E\sup_{t\in T}X'(t)\le C J_{ d}(D) \label{22a}
\end{equation}
for some $C<\infty$.
\end{remark}

\section{Local moduli of continuity}\label{sec-loc}

In this section we give a basic theorem for local moduli of continuity
of processes in $L^{\psi_{2}}$ in the spirit of Section~\ref{sec-2},
and apply it to permanental processes, as we do for the uniform modulus
of continuity in Section~\ref{sec-3}.

\begin{lemma}\label{lem-22q} Let $(T,\hat d)$ be a separable metric
or pseudometric space
with finite diameter $D$. Suppose that there exists
a probability measure $\mu$ on $(T,\hat d)$ such that $ J_{T,\hat
d,\mu} (D)<\infty$.

For any $t_{0}\in T$ and $\delta>0$, let $T_{\delta}:=\{s\dvtx \hat
d(s,t_{0})< \delta/2\}$. Suppose $0<\delta\le\delta_{0}<D$ which
implies that $T_{\delta}\subseteq T_{D}$. Consider the probability
measures $\mu_{\delta}(\cd):=\mu(\cd\cap T_{\delta})/\mu
(T_{\delta})$, $0<\delta\le\delta_{0}$, and assume that $c_{\mu
_{\delta},T_{\delta}}(\wt f)<\infty$, for each $0<\delta\le
\delta_{0}$; see~(\ref{c}) for the definition of $c_{\mu_{\delta
},T_{\delta}}$. Then
%
\begin{equation}\label{217w}\qquad
\sup_{\hat d(s,t_{0})< \delta/2} \bigl|f(s)-f(t_{0})\bigr|\le20
\sup_{t\in T_{\delta}}\int_{0}^{ \delta/4} \biggl(\log
\biggl(\frac{ c_{\mu_{\delta},T_{\delta}}(\wt f)+1}{\mu_{\delta}^{2}(B_{\hat d}(t,u))} \biggr)
 \biggr)^{1/2} \,du .
\end{equation}
\end{lemma}

\begin{pf}The condition that $ J_{T,\hat d,\mu} (D)<\infty$ implies
that $\mu(B_{\hat d}(t,u))>0$ for all $t\in T$ and $u>0$. Since
$T_{\delta}$ is open for every $t\in T_{\delta}$, there exists a
ball, say $B'_{\hat d}(t,u)\subset T_{\delta}$. Consequently,
%
\begin{equation}
\mu_{\delta}\bigl(B'_{\hat d}(t,u)\bigr)=
\frac{\mu( B'_{\hat d}(t,u))}{ \mu
(T_{\delta})}>0
\end{equation}
for all $t\in T$ and $u>0$. Therefore, (\ref{217w}) follows from
Lemma~\ref{lem-22}.
\end{pf}

The next corollary and theorem follow immediately from Lemma~\ref{lem-22q}.

\begin{corollary} Let
%
\begin{eqnarray}\label{220}
\qquad&& H_{T_{\delta},\hat d, \mu_{\delta} ,\delta} (\wt f)
\nonumber
\\[-8pt]
\\[-8pt]
\nonumber
&&\qquad =\delta\bigl(\log\bigl(c_{\mu_{\delta},T_{\delta}}(\wt f) +1\bigr)\bigr)^{1/2}+
\sup_{t\in T_{\delta}}\int_{0}^{ \delta/4} \biggl(\log
\biggl(\frac{ 1}{\mu_{\delta} (B_{\hat d}(t,u))} \biggr) \biggr)^{1/2} \,du.
\end{eqnarray}
Under the hypotheses of Lemma~\ref{lem-22q},
%
\begin{equation}
\lim_{\delta\to0} \sup_{ \hat d(s,t_{0})\le\delta
/2 }\frac
{ |f(s)-f(t_{0})|}{H_{T_{\delta},\hat d, \mu_{\delta} ,\delta}
(\wt f)}\le30\qquad \mbox{a.s.}
\label{21ww}
\end{equation}
\end{corollary}

\begin{theorem}\label{theo-23}Under the hypotheses of Theorem \ref
{maj}, assume
that (\ref{25}) holds. Define $\mu_{\delta}$ and $T_{\delta}$ as
in Lemma~\ref{lem-22q}. Then
%
\begin{equation}
\lim_{\delta\to0} \sup_{ \hat d(s,t_{0})\le\delta/2 }\frac{
|X'(s)-X'(t_{0})|}{H_{T_{\delta},\hat d, \mu_{\delta} ,\delta}
(\wt X)}\le30 \qquad\mbox{a.s.}
\label{21wwaa}
\end{equation}
\end{theorem}

We can use Theorem~\ref{theo-23} to find local moduli of continuity
for permanental processes. However, before we do this, we show that
with an additional mild regularity condition we can simplify the
expression in the denominator of (\ref{21wwaa}).
Consider the first term on the right-hand side of (\ref{220}), with
$\wt f$ replaced by $\wt X$.\vspace*{1pt} It is simply bounded by a constant times
$\delta$ unless $\limsup_{\delta\to0}c_{\mu_{\delta
},T_{\delta}}(\wt X)=\infty$ on a set of positive measure. Let us
assume this is the case. As in (\ref{213}), $ E c_{\mu_{\delta
},T_{\delta}}(\wt X) \leq1$. Therefore, for $\ep>0$,
%
\begin{eqnarray}
P \bigl(\log c_{\mu_{\delta},T_{\delta}}(\wt X) \ge(1+\ep )u \bigr)&\le& P \bigl(
c_{\mu_{\delta},T_{\delta}}(\wt X) \ge e^{(1+\ep) u } \bigr)
\nonumber
\\[-8pt]
\\[-8pt]
\nonumber
&\le&e^{- (1+\ep) u }.
\end{eqnarray}
It follows from the Borel--Cantelli lemma that for all $\bb<1$,
%
\begin{equation}
\limsup_{k\to\infty}\frac{\log c_{\mu_{\bb^{k}},T_{\bb
^{k}}}(\wt X)}{\log\log1/\bb^{k}}\le1.\label{47}
\end{equation}
We would like to extend this to get
%
\begin{equation}
\limsup_{\delta\to0}\frac{\delta(\log c_{\mu_{\delta
},T_{\delta}}(\wt X))^{1/2}}{\delta(\log\log1/\delta)^{1/2}}\le C.\label{225c}
\end{equation}
Note that for $ \bb^{k+1}<\delta\le\bb^{k}$,
%
\begin{eqnarray}
c_{\mu_{\delta},T_{\delta}} (\wt X)&=&\frac{1}{\mu^{2}(T_{\delta})}\int_{T_{\delta}\times T_{\delta}}
\psi_{2}(\wt X) \,d\mu(s ) \,d\mu(t)
\nonumber\\
&\leq&\frac{1}{\mu^{2}(T_{\delta})}\int_{T_{\bb^{k}}\times
T_{\bb^ {k}}}\psi_{2}(\wt X)
\,d\mu(s ) \,d\mu(t)
\\
&\le&\frac{\mu^{2}(T_{\bb^{k}})}{\mu^{2}(T_{\bb^{k+1}})} c_{\mu
_{\bb^{k}},T_{\bb^{k}}} (\wt X).
\nonumber
\end{eqnarray}
Consequently, if
%
\begin{equation}
\limsup_{k\to\infty}\frac{\mu(T_{\bb^{k}})}{\mu(T_{\bb^{k+1}}
)}\le C,\label{227}
\end{equation}
we can use (\ref{47}) to get (\ref{225c}).\vadjust{\goodbreak}

When (\ref{227}) holds we
have the following results for the local moduli of continuity of
permanental processes.

\begin{theorem}\label{theo-42} Under the hypotheses of Theorem \ref
{theo-11a}, assume that
(\ref{18jv}) and~(\ref{227}) hold. Then if $\th_{t_{0}}\ne0$
almost surely, there exists a version $\th'=\{\th'_{x},x\in T \}$
such that
%
\begin{equation}
\lim_{\delta\to0} \sup_{ d(s,t_{0})\le\delta/2 }\frac{ |\th'_{s}-\th'_{t_{0}}|}{\ov H_{T_{\delta},d,\mu_{\delta}} (\delta
/4)}\le C
\th^{1/2}_{t_{0}} \qquad\mbox{a.s.},\label{21wwdd}
\end{equation}
where
%
\begin{equation}\label{220aa}
 \ov H_{T_{\delta},d,\mu_{\delta}}(\delta/4) :=\delta(\log\log1/
\delta)^{1/2}+J_{T_{\delta},d,\mu_{\delta
}} (\delta/4).
\end{equation}
(See Lemma~\ref{lem-22q} for the definitions of the other terms.)

If $\th_{t_{0}}\equiv0$, there exists a version $\th'=\{\th'_{x},x\in T \}$ such that
%
\begin{equation}
\lim_{\delta\to0} \sup_{ \hat d(s,t_{0})\le\delta/2 }\frac{
\th'_{s} }{ (\ov H_{T_{\delta},d,\mu_{\delta}} (\delta/4)
)^{2}}\le C\qquad \mbox{a.s.}
\label{21qqq}
\end{equation}
\end{theorem}

\begin{pf}We use Theorem~\ref{theo-23} with $X=\th^{1/2}$ and (\ref
{225c}) and the same argument used in the proof of Theorem \ref
{theo-11a}, in particular (\ref{ffun}), to get
%
\begin{equation}
\lim_{\delta\to0} \sup_{ d(s,t_{0})\le\delta/2 }\frac{ |\theta^{1/2}_{s}-\theta^{1/2}_{t_{0}}|}{\ov H_{T_{\delta},d,\mu_{\delta
}} (\delta/4)}\le C \qquad\mbox{a.s.}
\label{21swa}
\end{equation}
It is easy to see that this gives (\ref{21qqq}). To get (\ref{21wwdd}),
fix $\delta'>0$. Then for any $\delta\leq\delta'$
\[
\sup_{ d(s,t_{0})\le\delta/2 }|\theta_{s}-\theta_{t_{0}}|\leq \bigl|
\theta^{1/2}_{s}-\theta^{1/2}_{t_{0}}\bigr| 2
\sup_{z\in
B_{d}(t_{0},\delta')}\theta^{1/2}_{z},
\]
so we obtain
%
\begin{equation}
\lim_{\delta\to0} \sup_{ d(s,t_{0})\le\delta/2 }\frac{ |\theta_{s}-\theta_{t_{0}}|}{\ov H_{T_{\delta},d,\mu_{\delta}} (\delta
/4)}\le C
\sup_{z\in B_{d}(t_{0},\delta')}\theta^{1/2}_{z} \qquad\mbox{a.s.}\label{21swh}
\end{equation}
Letting $\delta'\rar0$
completes the proof.
\end{pf}

\begin{remark} Note that if $\th$ is the square of Gaussian
process, $\ov H_{T_{\delta},d,\mu_{\delta}} (\cd) $ is equivalent
to the correct local modulus of continuity of the Gaussian process.
\end{remark}

\begin{example}\label{ex-21} Theorems~\ref{theo-23} and \ref
{theo-42} seem very abstract. We show here how they give the familiar
iterated logarithm behavior for fairly regular processes on nice spaces.

Take $T$ to be the unit interval in $R^{1}$. Assume that
%
\begin{equation}
\hat d(s,t_{0})= \phi\bigl(|s-t_{0}|\bigr) \qquad\mbox{for }
0<|s-t_{0}|\le\delta_{0}\label{227s}\vadjust{\goodbreak}
\end{equation}
for some $\delta_{0}>0$, and some continuous increasing function
$\phi$. Now take $\mu$ to be Lebesgue measure. In this case,
%
\begin{equation}
\mu(T_{\delta})= 2\phi^{-1}(\delta/2),\label{228}
\end{equation}
so that, for example, (\ref{227}) holds if $\phi$ is regularly
varying. In addition, it follows from~\cite{book}, (7.94), that the
second term on the right-hand side of (\ref{220}), with $\wt f$
replaced by $\wt X$, is bounded by a constant times
%
\begin{equation}
\delta+\int_{0}^{1}\frac{\phi(\phi^{-1}(\delta/2)u)}{u(\log
2/u)^{1/2}} \,du.
\end{equation}
Note that under (\ref{227s}) we can replace $\hat d(s,t_{0})\le
\delta/2$ in (\ref{21wwaa}) by $|s-t_{0}|\le\phi^{-1}(\delta
/2)$. Then, replacing $ \phi^{-1}(\delta/2)$ by
$\delta'$ and making a change of variables, as in~\cite{book},
(7.96), and using (\ref{225c}), we get
%
\begin{equation}
\lim_{\delta'\to0} \sup_{ |s-t_{0}|\le\delta' }\frac{
|X'(s)-X'(t_{0})|}{\wt H(\delta') }\le C \qquad\mbox{a.s.},
\label{21wwz}
\end{equation}
where
%
\begin{equation}
\wt H(\delta)=\phi(\delta) \bigl(\log\log1/\phi(\delta/2)\bigr)^{1/2} +
\int_{0}^{1}\frac{\phi(2\delta u)}{u(\log2/u)^{1/2}} \,du.
\end{equation}
By~\cite{book}, (7.128), if $\phi$ is regularly varying,
%
\begin{equation}
\lim_{\delta\to0} \frac{\wt H(\delta)}{\phi(\delta)(\log\log
1/ \delta)^{1/2}}=1.
\end{equation}

In the same vein, under (\ref{227}) and the assumption that $\phi$
is regularly varying, it follows from (\ref{theo-22}) and the
material in~\cite{book}, pages 298 and 299, that
%
\begin{equation}
\lim_{\delta\to0} \sup_{ |s-t |\le\delta}\frac{ |X'(s)-X'(t
)|}{\phi(\delta)( \log1/ \delta)^{1/2} }\le C \qquad\mbox{a.s.}
\label{21wwzf}
\end{equation}
\end{example}

\section{Dominating metrics for permanental processes}\label{sec-4}

We exhibit several interesting metrics and other functions that
dominate $d$ or are even equivalent to $d$. [$ d_{1}$ is equivalent to
$d$ ($ d\approx d_{1}$) if there exist constants $0<c_{1}\le
c_{2}<\infty$ such that $c_{1}\,d\le d_{1}\le c_{2}\,d$.] Note that for
$C\ne0$,
%
\begin{equation}
J_{T,Cd,\mu}(a)=C J_{T, d,\mu}(a/C).\label{41}
\end{equation}
Therefore, multiplying a metric or related function by a constant
alters our results in an acceptable way.

We consider several scenarios. To simplify the exposition we work with
%
\begin{eqnarray}\label{52}
\ov d(x,y)&:=& d(x,y) /4\sqrt{ 2/3}
\nonumber
\\[-8pt]
\\[-8pt]
\nonumber
&\hphantom{:}= &\bigl(\Ga(x,x)+\Ga(y,y)-2 \bigl(\Ga(x , y)
\Ga(y, x ) \bigr)^{1/2} \bigr)^{1/2}.
\end{eqnarray}

(1) \textit{Conditions under which ${\ov{d}}$ is equivalent to
natural metrics for $\th$}.

\begin{lemma}\label{lem-d2} Let
%
\begin{equation}
d_\th(x,y)= \bigl( E ( \th_{x} - \th_{y}
)^{2} \bigr)^{1/2}
\end{equation}
and
%
\begin{equation}
\hat d_\th(x,y)= \bigl( E \bigl((\th_{x}-E
\th_{x})-(\th_{y}-E\th_{y}) \bigr)^{2}
\bigr)^{1/2}.
\end{equation}
Then
%
\begin{equation}
\tfrac{\sqrt2}{\sqrt2+1}\,d_\th(x,y)\le\hat d_\th(x,y)\le2
\,d_\th (x,y)\label{wqw}
\end{equation}
and
%
\begin{eqnarray}\label{met4}
K \bigl(\Ga(x,x)+\Ga(y,y) \bigr)^{1/2} \ov d (x,y)&\leq&\hat
d_\th (x,y)
\nonumber
\\[-8pt]
\\[-8pt]
\nonumber
&\leq& 2 \bigl(\Ga(x,x)+\Ga(y,y) \bigr)^{1/2} \ov d
(x,y),
\end{eqnarray}
where $K=\sqrt2/(\sqrt2+1)$.
\end{lemma}

\begin{remark} By (\ref{met4}),
%
\begin{equation}
c_{1} \ov d (x,y)\leq\hat d_\th(x,y)\leq c_{2}
\ov d (x,y),\label{met4a}
\end{equation}
where $c_{1}=2\inf_{x\in T}\Ga^{1/2}(x,x),  c_{2}=2\sqrt2 \sup_{x\in T}\Ga^{1/2}(x,x) $.
In particular, if
\[
0<\inf_{x\in T}\Ga(x,x) \leq\sup_{x\in T}\Ga(x,x)<\infty,
\]
then $d$ is equivalent to
$\hat d_\th$ and $d_\th$.
\end{remark}

\begin{pf}
By Lemma~\ref{lem-1},
%
\begin{equation}\label{48w}
\hat d_\th^{2}(x,y)=2 \bigl(\Ga^{2}(x,x)+
\Ga^{2}(y,y ) -2\Ga(x,y) \Ga(y,x) \bigr) .
\end{equation}
Let
%
\begin{equation}
\qquad\wt d^{2}(x,y):= (E\th_{x}-E\th_{y}
)^{2}= \bigl(\Ga^{2}(x,x)+ \Ga^{2}(y,y ) -2
\Ga(x,x) \Ga(y,y) \bigr).\label{ka11}
\end{equation}
By (\ref{32}),
%
\begin{equation}
\wt d (x,y)\le\tfrac{1}{\sqrt{2}} \hat d_\th(x,y).\label{515m}
\end{equation}
By the Cauchy--Schwarz inequality,
%
\begin{equation}
\wt d (x,y)\le d_\th(x,y).
\end{equation}
Using this and the triangle inequality, we see that
%
\begin{equation}
\hat d_\th(x,y)\le d_\th(x,y)+\wt d(x,y)\le2
\,d_\th(x,y)
\end{equation}
and
%
\begin{equation}
\hat d_\th(x,y)\ge d_\th(x,y)-\wt d(x,y),
\end{equation}
which, along with (\ref{515m}), implies that
%
\begin{equation}
\bigl(1+\tfrac{1}{\sqrt{2}} \bigr)\hat d_\th(x,y)\ge
d_\th(x,y).
\end{equation}
Thus we get (\ref{wqw}).

By (\ref{48w}) and (\ref{15aa}),
%
\begin{eqnarray}\label{48q}
\hat d_\th^{2}(x,y)&\le&2 \bigl(\bigl(\Ga(x,x)+ \Ga(y,y )
\bigr)^{2} -4\Ga(x,y) \Ga(y,x) \bigr)
\nonumber\\
&=& 2 \bigl(\Ga(x, x)+\Ga(y, y)-2\sqrt{\Ga(x, y)\Ga(y, x)} \bigr)
\\
&&{} \times \bigl(\Ga(x, x)+\Ga(y, y)+2\sqrt{\Ga(x, y)\Ga(y, x)} \bigr).
\nonumber
\end{eqnarray}
This gives the upper bound in (\ref{met4}).

For the lower bound, we note that
%
\begin{eqnarray}\label{met1}
d_\th^{2}(x,y) &=&E \bigl(G^{2}(x )
-G^{2}(y) \bigr)^{2}
\nonumber
\\
&=& E \bigl\{ \bigl(G(x ) -G(y) \bigr)^{2} \bigl(G(x )+G(y)
\bigr)^{2} \bigr\}
\nonumber\\
&=& E \bigl(G(x ) -G(y) \bigr)^{2} E \bigl(G(x )+G(y)
\bigr)^{2} +2 \bigl( E \bigl\{ G^{2}(x )
-G^{2}(y) \bigr\} \bigr)^{2}
\nonumber
\\[-8pt]
\\[-8pt]
\nonumber
&\ge&E \bigl(G(x ) -G(y) \bigr)^{2} E \bigl(G(x )+G(y)
\bigr)^{2}
\\
&=& \bigl(\Ga(x, x)+\Ga(y, y)-2\sqrt{\Ga(x, y)\Ga(y, x)} \bigr)
\nonumber
\\
&&{} \times \bigl(\Ga(x, x)+\Ga(y, y)+2\sqrt{\Ga(x, y)\Ga(y, x)} \bigr)
\nonumber
\end{eqnarray}
Consequently,
%
\begin{equation}
d_\th(x,y) \ge\bigl(\Ga(x, x)+\Ga(y, y)\bigr)^{1/2} \ov
d(x,y).
\end{equation}
Using (\ref{wqw}) we get the lower bound in (\ref{met4}).
\end{pf}

\begin{lemma}\label{lem-metlb}
%
\begin{equation}
\ov d(x,y)\le d_\th^{1/2}(x,y)\label{520}.
\end{equation}
\end{lemma}

\begin{pf} By (\ref{met1}),
%
\begin{equation}
d_\th^{2}(x,y) \ge \bigl(\bigl(\Ga(x, x)+\Ga(y, y)
\bigr)^{2}-4 \Ga(x, y)\Ga(y, x) \bigr) .
\end{equation}
Consequently,
%
\begin{equation}
d_\th(x,y) \ge \bigl(\bigl(\Ga(x, x)+\Ga(y, y)\bigr) -2 \bigl(
\Ga(x, y)\Ga(y, x) \bigr)^{1/2} \bigr) .
\end{equation}
Taking the square root again, we get (\ref{520}).
\end{pf}

\begin{lemma}\label{lem-cmet}
%
\begin{eqnarray}\quad
\bigl| \ov d(x,y)-\ov d(x,z)\bigr|\le C \Bigl(1+\sup_{u\in T}\Ga(u,u) \Bigr) \bigl(
\hat d_\th^{1/4}(y,z) +\hat d_\th^{1/2}(y,z)
\bigr).\label{cmet}
\end{eqnarray}
\end{lemma}

\begin{pf}
%
\begin{eqnarray}
&&\bigl| \ov d(x,y)-\ov d(x,z)\bigr|
\nonumber
\\
&&\qquad \le\bigl| \ov d^{2}(x,y)-\ov d^{2}(x,z)\bigr|^{1/2}
\nonumber
\\[-8pt]
\\[-8pt]
\nonumber
&&\qquad \le\bigl|\Ga(y,y)-\Ga(z,z)\bigr|^{1/2}+2\bigl|\bigl(\Ga(x,z)\Ga (z,x)
\bigr)^{1/2}-\bigl(\Ga(x,y)\Ga(y,x)\bigr)^{1/2}\bigr|^{1/2}
\\
&&\qquad \le\bigl|\Ga(y,y)-\Ga(z,z)\bigr|^{1/2}+ 2\bigl| \Ga(x,z)\Ga(z,x) - \Ga (x,y)
\Ga(y,x) \bigr|^{1/4}
.\nonumber
\end{eqnarray}
By
(\ref{48w}),
%
\begin{eqnarray}
&& \bigl| \Ga(x,z)\Ga(z,x) - \Ga(x,y)\Ga(y,x)\bigr |
\nonumber
\\[-8pt]
\\[-8pt]
\nonumber
&&\qquad \le C \bigl(\bigl|\hat d_\th^{2}(x,z)-\hat
d_\th^{2}(x,y)\bigr|+\bigl|\Ga^{2}(y,y)-
\Ga^{2}(z,z)\bigr| \bigr),
\end{eqnarray}
and by (\ref{515m}),
%
\begin{equation}
\bigl|\Ga(y,y)-\Ga(z,z)\bigr|\le\hat d_\th(y,z).
\end{equation}
In addition,
%
\begin{eqnarray}
\bigl|\hat d_\th^{2}(x,z)-\hat d_\th^{2}(x,y)\bigr|&
\le& 2\sup_{u,v\in T}\hat d_\th(u,v) \bigl|\hat d_\th(x,z)-
\hat d_\th(x,y)\bigr|
\nonumber
\\[-8pt]
\\[-8pt]
\nonumber
&\le&8\sup_{u\in T}\Ga(u,u)\hat d_\th(y,z)
.
\end{eqnarray}
Putting these together we get (\ref{cmet}).
\end{pf}

\begin{lemma}\label{lem-cmetop} Assume that $\sup_{u\in T}\Ga
(u,u)<\infty$.
Then the sets $b_{\ov d}(x,u)=\{y\in T | \ov d(x,y)<u\}, x\in T,u\in
R_{+} $ form the base for the $\hat d_\th$ (and equivalently the
$d_{\th}$) metric topology.
\end{lemma}

\begin{pf} Let $f_{x}(y)=\ov d(x,y).$ By (\ref{cmet}) we have that $f_{x}$
is continuous with respect to
$\hat d_\th$, and hence $b_{\ov d}(x,u)=f^{-1}_{x}([0,u))$ is open
with respect to
$\hat d_\th$. We now show that for any $x\in T,u\in R_{+} $, and any
$y\in b_{\hat d_\th}(x,u)$, we can find $v>0$ such that $b_{\ov
d}(y,v)\subseteq b_{\hat d_\th}(x,u)$. To see this, first choose
$w>0$ such that $b_{\hat d_\th}(y,w)\subseteq b_{\hat d_\th}(x,u)$.
It then follows from (\ref{met4a}) that
$b_{\ov d}(y,c^{-1}_{2}w)\subseteq b_{\hat d_\th}(y,w)$. By (\ref
{wqw}) the same argument applies with $\hat d_\th$ replaced by $d_{\th
}$.
\end{pf}

Let $\Sigma(x,y)=\Ga(x,y)\Ga(y,x)$. It follows from Lemma \ref
{lem-1} that $\{\Sigma(x,y), x,y\in T\}$ is positive definite.
Therefore it is the covariance of a mean zero Gaussian process which we
denote by $\{\SS(x),x\in T\}$. Clearly,
%
\begin{equation}
\hat d_\th(x,y)= \bigl( E (\SS_{x}-\SS_{y}
)^{2} \bigr)^{1/2}.
\end{equation}

(2) \textit{Conditions under which ${\ov{d}}$ is equivalent to a
function that may be a metric for a Gaussian process}. We suppose that
%
\begin{equation}
\bigl| \Ga(x, y)\bigr|\vee\bigl|\Ga(y,x)\bigr| \le\Ga(y,y) \wedge\Ga(x, x) .\label{49}\vadjust{\goodbreak}
\end{equation}
Let
%
\begin{equation}
d_{2}(x,y)= \bigl\{ \Ga(x, x)+\Ga(y, y)-\bigl(\bigl|\Ga(x, y)\bigr|+\bigl|\Ga(y, x)
\bigr|\bigr) \bigr\}^{1/2}.\label{410}
\end{equation}

\begin{lemma}\label{lem-53} When (\ref{49}) holds,
%
\begin{equation}
\tfrac{1}{\sqrt2} \ov d (x,y)\le d_{2}(x,y)\le\ov d (x,y).
\label{411}
\end{equation}
\end{lemma}

In general when $\Ga(x,y)$ is the potential density of a Borel right
process $X$, in place of (\ref{49}), we only have
%
\begin{equation}
0\leq\Ga(x, y)\le\Ga(y,y)\quad \mbox{and}\quad 0\leq\Ga (y,x) \le\Ga(x, x); \label{49a}
\end{equation}
see, for example,~\cite{book}, Lemma 3.3.6, where this is proved for
symmetric potential densities, and note that the proof also works when
the densities are not symmetric.

Set
$\wt\Ga(x,y)=\Ga(y,x)$. This is the potential density of $\wt X$,
the dual process of $X$. Therefore, if $\wt X$ is also a Borel right
process, using (\ref{49a}), we actually get (\ref{49}). In \cite
{EK1} it is shown that for certain Borel right processes $X$ with
potential density $\Ga(x,y)$, the symmetric function $\frac{\Ga
(x,y)+\Ga(y,x)}{2}$ is positive definite, so that $d_{2}(x,y)$ is the
$L^{2}$ metric of a Gaussian process; see Section~\ref{sec-comp} for details.

\begin{pf*}{Proof of Lemma~\ref{lem-53}} We have
%
\begin{eqnarray}\label{412}
\ov d^ { 2} (x,y)&=&d_{2}^ { 2} (x,y)+\bigl| \bigl|
\Ga(x,y)\bigr|^{1/2}- \bigl|\Ga (y,x)\bigr|^{1/2}\bigr|^{2}
\nonumber
\\[-8pt]
\\[-8pt]
\nonumber
&\le& d_{2}^ { 2} (x,y)+\bigl| \bigl|\Ga(x,y)\bigr| - \bigl|\Ga(y,x)\bigr| \bigr| .
\end{eqnarray}
By (\ref{49}) if $ |\Ga(x,y)| - |\Ga(y,x)|\ge0$, then
%
\begin{eqnarray}
\bigl|\bigl|\Ga(x,y)\bigr| - \bigl|\Ga(y,x)\bigr|\bigr|&\le&\Ga(y,y) - \bigl|\Ga(y,x)\bigr|
\nonumber\\
&\le&
\Ga(y,y) +\Ga(x,x)- \bigl(\bigl|\Ga(y,x)\bigr| +\bigl|\Ga (x,y)\bigr| \bigr)
\\
&=&d_{2}^ { 2} (x,y).
\nonumber
\end{eqnarray}
Interchanging $x$ and $y$, we also get that when
and if $|\Ga(y,x)|- |\Ga(x,y)|\ge0$,
%
\begin{equation}
\bigl|\Ga(y,x)\bigl|- \bigl|\Ga(x,y)\bigr| \le d_{2}^ { 2} (x,y).
\end{equation}
Therefore,
%
\begin{equation}
\ov d^ { 2} (x,y)\le2 \,d_{2}^ { 2} (x,y).
\end{equation}
Using this and the first line of (\ref{412}), we get (\ref
{411}).
\end{pf*}

\section{Local times of Borel right processes}

Our primary motivation for obtaining sample path properties of
permanental processes was to use them, along with the following
isomorphism theorem, to obtain sample path properties\vadjust{\goodbreak} of the local
times of Borel right processes, paralleling our use of Dynkin's
isomorphism theorem in~\cite{sip}, to obtain sample path properties of
the local times of strongly symmetric Borel right processes.

Let $X =
(\Om, X_t ,P^x
)$ be a Borel right process with
$0$-potential density
$u(x,y)$. Let $h_{x}(z)=u(z,x)$, and assume that $h_{x}(z)>0$ for all
$x,z\in S$. Recall that the expectation operator $E^{z/h_{x}}$ for the
$h_{x}$-transform of $X$ is given by
%
\begin{equation}
E^{z/h_{x}}(F 1_{\{t<\zeta\}})=\frac1{h_{x}(z)}E^z
\bigl(F h_{x}(X_t )\bigr) \label{sp2}
\end{equation}
for all bounded $\mathcal{F}^{0}_t$ measurable functions $F$,
where $\mathcal{F}^{0}_t$ is the $\sigma$-algebra generated by $\{
X_{r}, 0\leq r\leq t\}$; see, e.g.,~\cite{book}, (3.211). Here, as
usual, $E^z$ denotes the expectation operator for $X$ started at $z$.

Recall that on page \pageref{page-3} we wrote that Eisenbaum
and Kaspi pointed out that the $0$-potential of a transient Markov
process was a kernel for a permanental process. Using this they
establish the following isomorphism theorem.

\begin{theorem}[(Eisenbaum and Kaspi~\cite{EK})]\label{theo-ke} Let $X
=
(\Om, X_t ,P^x
)$ be a Borel right process with
$0$-potential density
$u(x,y)$, and let
$L=\{L^y_{t} ; (y,t)\in S\times R_+\}$
denote the local times for
$X$, normalized so that
%
\begin{equation}
E^{ v} \bigl( L^{y}_\infty \bigr)=u(v,y).
\label{mp485}
\end{equation}
Let $x$ denote a fixed element of $S$, and assume that $u(x,x)>0$. Set
%
\begin{equation}
h_{x}(z)= u(z,x).\label{sp1}
\end{equation}
Let $\th=\{\th_{y} ; y\in S\}$ denote the
permanental process with kernel $u(x,y)$. Then, for any
countable subset $D\subseteq S$,
%
\begin{equation}
\biggl\{ L^y_{\infty}+\frac{ 1}{2}\th_{y} ; y
\in D , P^{x/h_{x}}\times P_{\th} \biggr\}\stackrel{\mathrm{law}} {=} \biggl
\{\frac{ 1}{2}\th_{y} ; y\in D , \frac{\th_{x}}{u(
x,x)}P_{\th}
\biggr\}.\label{it18}
\end{equation}

Equivalently, for all $x_{ 1},\ldots,x_{ n}$ in $S$ and bounded
measurable functions
$F$ on $R^n_+$, for all
$n$,
%
\begin{equation}
E^{x/h_{x}}E_{ \th} \biggl(F \biggl(L^{x_{ i}}_{\infty}+
\frac{ \th_{x_{ i}} }2 \biggr) \biggr) = E_{ \th} \biggl(\frac{\th_{x}}{u(
x,x)}F
\biggl( \frac{ \th_{x_{ i}} }2 \biggr) \biggr).\label{it19}
\end{equation}
[Here we use the notation $F( f(x_{ i})):=F( f(x_{ 1}),\ldots,
f(x_{ n}))$.]
\end{theorem}

Theorem~\ref{theo-ke} is only a partial analog of Dynkin's isomorphism
theorem for strongly symmetric Borel right processes,~\cite{book}, Theorem
8.1.3, which holds with measures $P^{x/h}$, for a much wider
class of functions $h$ than those in (\ref{sp1}).
In addition, note that
Theorem~\ref{theo-ke} can only give a version of
$\{L^y_{t} ; (y,t)\in S\times R_+\}$ which is jointly continuous with
respect to the measures $P^{x/h_{x}}$. In order to use this to obtain
joint continuity with respect to the measures $P^{x}$, we use (\ref
{sp2}) with $z=x$.
Therefore, since we require that $h_{x}(z)>0$ for all $z\in S$, when
$P^{x/h_{x}}(A, t<\zeta)=0$ for some $A\in\mathcal{F}_t^{0}$,
we also have $P^{x }(A, t<\zeta)=0$.\vadjust{\goodbreak}

When we say that a stochastic process
$\hat L=\{\hat L^y_t,(y,t)\in S \times R_+ \}$
is a version of the local time of
a Markov process $X$ we mean more than the traditional statement that
one stochastic process is a version of the other. Besides this, we also
require that the version is itself a local time for $X$, that is, that for
each $y\in S$, $\hat L^y_\cdot$ is a local time for $X$ at $y$. To be
more specific, suppose that
$L=\{L^y_t,(y,t)\in S \times R_+\}$ is a local time for $X$. When we say
that we can find a version of the local time which is jointly continuous
on $S\times T$, where $T\subset R_+$, we mean that
we can find a stochastic process
$\hat L=\{\hat L^y_t,(t,y)\in(y,t)\in S \times R_+\}$
which is continuous on $ S\times T$ for all $x\in
S$ and which satisfies, for each $x,y\in S$
%
\begin{equation}
\hat L^y_t=L^y_t\qquad  \forall t\in
R_+, P^x \mbox{ a.s. }
\end{equation}
Following convention, we often say that a Markov process
has a continuous local time, when we mean that we can find a continuous
version for the local time.

\begin{pf*}{Proof of Theorem~\ref{theo-jc}}
The proof
follows the general lines of the proof for symmetric Markov processes
in~\cite{sip}, Section 6. However, there are significant differences,
so we give a self-contained proof.

Since $S$ is a locally compact topological space with a countable base,
we can find a metric $\rho$ which induces the topology of $S$.
We first consider the case where $X$ is a transient Borel right process
with state space $S$ and continuous, strictly positive
$0$-potential densities
$u (x,y)$. We take $\th$ to be the permanental process with kernel $u (x,y)$.

Fix a compact set $K\subseteq T$ and some $x\in K$. By (\ref{18jj}),
Theorems~\ref{theo-11a} and~\ref{maj} we can find a version
of $\th$ which is continuous on $K$ almost surely, and such that for
each $p$,
%
\begin{equation}
E\sup_{x\in K}\th^{ p}_{x} <\infty. \label{momentp}
\end{equation}
We work with this version.

It follows
from~\cite{sip}, (4.30) and (4.31), that for any $z,y\in S$
%
\begin{equation}
E^{z/h_{x}}\bigl(L^y_\infty\bigr) =\frac{u(z,y)h_{x}(y)}{h_{x}(z)}.
\label{sip2c}
\end{equation}
We shall use the fact that that $X_{t}$ is a right continuous simple
Markov process under the measures $P^{z/h_{x}}$,~\cite{book}, Lemma 3.9.1.

To begin, we first show first that $L$ is jointly continuous on
$K\times R_+ $, almost
surely with respect to $P^{x/h_{x}}$. By~\cite{book}, Lemma 3.9.1, we
can assume that the local times $L^y_t$ are $\mathcal F^{0}_t$ measurable.
Consider the martingale
%
\begin{equation}
A^y_t=E^{x/h_{x}}\bigl(L^y_{\infty}
\mid\mathcal F^{0}_t\bigr) .
\end{equation}
Let $\tau_{t}$ denote the shift operator on $\Om$. Then
%
\begin{equation}
L^y_{\infty}=L^y_t+L^y_{\infty}
\circ\tau_t =L^y_t+1_{\{t<\ze\}}L^y_{\infty}
\circ\tau_t.\vadjust{\goodbreak}
\end{equation}
Therefore
%
\begin{eqnarray}
A^y_t&=&L^y_t+E^{x/ h_{x}}
\bigl(1_{\{t<\ze\}}L^y_{\infty}\circ\tau_t\mid
\mathcal F^{0}_t\bigr)
\nonumber
\\[-8pt]
\\[-8pt]
\nonumber
&=& L^y_t+1_{\{t<\ze\}} E^{x/ h_{x}}
\bigl(L^y_{\infty}\circ\tau_t\mid \mathcal
F^{0}_t\bigr) =L^y_t+1_{\{t<\ze\}}
E^{X_t/ h_{x}}\bigl(L^y_{\infty}\bigr),
\end{eqnarray}
where we use the simple Markov property described above.
It follows from (\ref{sip2c}), using
the convention that $1/h(\De)=0$, that
%
\begin{equation}
A^y_t=L^y_t+
\frac{u(X_t,y)h_{x}(y)}{h_{x}(X_t)}. \label{1200}
\end{equation}
Since $X_t$ is right continuous for $P^{x/h_{x}}$,
$A_t^y$ is also right continuous.
Let $D$ be a countable, dense subset of $K$, and $F$ a finite
subset of $D$. Since
%
\begin{equation}
\mathop{\sup_{\rho(y,z)\le\delta}}_{y,z\in F} A^y_t-A^z_t
=\mathop{\sup_{\rho(y,z)\le\delta}}_{y,z\in F} \bigl|A^y_t-A^z_t\bigr|
\end{equation}
is a right continuous, nonnegative submartingale, we have, for any
$\ep>0$,
%
\begin{eqnarray}\label{120}
&& P^{x/h_x}\Bigl(\sup_{t\ge0}\mathop{\sup_{\rho(y,z)\le\delta}}_{y,z\in F}
A^y_t-A^z_t\ge\ep\Bigr)
\nonumber
\\[-8pt]
\\[-8pt]
\nonumber
&&\qquad \le\frac1{\ep}E^{x/h_x}\Bigl( \mathop{\sup_{\rho(y,z)\le
\delta}}_{y,z\in F}
L^y_{\infty}-L^z_{\infty}\Bigr) \le\frac1{
\ep}E^{x/h_x}\Bigl(\mathop{\sup_{\rho(y,z)\le\delta}}_
{y,z\in D} L^y_{\infty}-L^z_{\infty}
\Bigr).
\end{eqnarray}
It
follows from (\ref{it19}) that
%
\begin{eqnarray}\label{121}
E^{x/h_x} \Bigl(\mathop{\sup_{\rho(y,z)\le\delta}}_{y,z\in D} L^y_{\infty}-L^z_{\infty}
\Bigr) &\le & E_\th \biggl(\mathop{\sup_{\rho(y,z)\le\delta}}_{y,z\in D} \biggl|\frac{\th_{y}}2-
\frac{\th_{z}}2 \biggr| \biggr)
\nonumber
\\[-8pt]
\\[-8pt]
\nonumber
&&{} +\frac{1 }{u(x,x)} \biggl(E_\th \biggl(\mathop{\sup_{\rho
(y,z)\le\delta}}_{y,z\in D} \biggl|
\frac{\th_{y}}2-\frac{\th_{z}}2 \biggr|^2 \biggr)E_\th
\bigl(\th^{2}_{x}\bigr) \biggr)^{1/2}.
\end{eqnarray}

It
follows from the uniform continuity of $\th$ on $K$ and (\ref
{momentp}) that for any $\bar\ep>0$, we can choose a
$\delta>0$ such that the right-hand side (\ref{121}) is less
than $\bar\ep$.
Combining~(\ref{1200})--(\ref{121}), we get
%
\begin{eqnarray} \label{122a}
&& P^{x/h_x}\biggl(\sup_{t\ge0}\mathop{\sup_{\rho(y,z)\le\delta}}_{y,z\in F}
L^y_t-L^z_t\ge2\ep
\nonumber\hspace*{-15pt}
\\
&&\hspace*{-20pt}\qquad\qquad \le\bar\ep+P^{x/h_x} \biggl(\sup_{t\ge0}\frac
1{h(X_t)}\mathop{\sup_{\rho(y,z)\le\delta} }_{y,z\in D}  \bigl(u(X_t,y)h_{x}(y)-u(X_t,z)h_{x}(z)
\bigr)\ge\ep \biggr)\biggr)\hspace*{-15pt}
\\
&&\qquad \le\bar\ep+P^{x/h_x} \biggl(\sup_{t\ge0}\frac1{h_{x}(X_t)}
\ge \frac{\ep}{\gamma(\delta)} \biggr),
\nonumber\hspace*{-15pt}
\end{eqnarray}
where
%
\begin{eqnarray}\label{122b}
\gamma(\delta)&=&\sup_{w\in S}\mathop{\sup_{\rho(y,z)\le
\delta}}_{y,z\in D}\bigl |u(w,y)h_{x}(y)-u(w,z)h_{x}(z)\bigr|
\nonumber
\\[-8pt]
\\[-8pt]
\nonumber
&=&\sup_{w\in K}\mathop{\sup_{\rho(y,z)\le\delta}}_{y,z\in D} \bigl|u(w,y)h_{x}(y)-u(w,z)h_{x}(z)\bigr|.
\end{eqnarray}
The last equality follows from~\cite{book}, (3.69), since the proof
does not require that $u(x,y)$ is symmetric.

It follows easily from (\ref{sp2}) and the fact that $X_{t}$ is a
simple Markov process under the measures $P^{z/h_{x}}$, that
$1/h_{x}(X_t)$ is a
supermartingale with respect to $P^{x/h_x}$.
Since $1/h_{x}(X_t)$ is also right continuous and nonnegative, we have
%
\begin{eqnarray}\label{122}
 P^{x/h_x} \biggl(\sup_{t\ge0}\frac1{h_{x}(X_t)}
\ge \frac{\ep
}{\gamma(\delta)} \biggr)&\le& \frac{\gamma(\delta)}{\ep} E^{x/h_x} \biggl(\frac
1{h_{x}(X_0)} \biggr)=\frac{\gamma(\delta)}{\ep h_{x}(x)}
\nonumber
\\[-8pt]
\\[-8pt]
\nonumber
&=&
\frac
{\gamma(\delta)}{\ep}.
\end{eqnarray}
Since both $h$ and $u$ are bounded and uniformly continuous on $K$, it
follows from~(\ref{122b})
that by
choosing
$\delta>0$ sufficiently small, we can make the right-hand side of
(\ref{122})
less than $\bar\ep$. By this observation and (\ref{122a}), and
taking the
limit over a sequence of finite sets increasing to $D$, we see that for
any $\ep$ and $\bar\ep>0$, we can find a $\delta>0$ such that
\[
P^{x/h_x}\Bigl(\sup_{t\ge0}\mathop{\sup_{\rho(y,z)\le\delta}}_{y,z\in D}
L^y_t-L^z_t\ge2\ep\Bigr)\le2
\bar\ep.
\]
It follows by the Borel--Cantelli lemma that we can find a sequence
$\{\delta_i\}_{i=1}^{\infty}$, $\delta_i>0$, such that $\lim_{i\to\infty}\delta_i=0$
and
%
\begin{equation}\label{123}
\sup_{t\ge0}\mathop{\sup_{\rho(y,z)\le\delta_{i}}}_{y,z\in D} L^y_t-L^z_t
\le\frac1{2^i}
\end{equation}
for all $i\ge I(\omega)$, almost surely with respect to $P^{x/h_x}$.

Fix $T<\infty$. We will now show that $L^y_t$
is uniformly continuous on $[0,T]\times D$, almost surely, with
respect to $P^{x/h_x}$. That is, for
each $\omega\in\Om'\subseteq\Om$, with $P^{x/h_x}(\Om')=1$, we
can find an
$I(\omega)$, such that for $i\ge I(\omega)$,
%
\begin{equation} \label{124}
\mathop{\sup_{ |s-t|\le\delta'_i}}_{ s,t\in[0,T]} \mathop{\sup_{\rho(y,z)\le\delta'_i}}_{ y,z\in D } \bigl|L^y_s-L^z_t\bigr|
\le\frac1{2^i},
\end{equation}
where $\{\delta'_i\}_{i=1}^{\infty}$ is a sequence of real numbers
such that
$\delta'_i>0$ and \mbox{$\lim_{i\to\infty}\delta'_i=0$}.\vadjust{\goodbreak}

To prove (\ref{124}), fix $\omega$ and assume that $i\ge I(\omega
)$, so that (\ref{123})
holds. Let $Y=\{y_1,\ldots,y_n\}$ be a finite subset of $ D$ such
that
\[
K\subseteq\bigcup_{j=1}^nB_{\rho}(y_j,
\delta_{i+2}).
\]
By definition,
each $L^{y_j}_t(\omega)$, $j=1,\ldots,n$, is uniformly continuous on $[0,T]$.
Therefore
we can find a finite, increasing sequence
$t_1=0,t_2,\ldots,t_{k-1}<T,t_k\ge T$ such that
$t_m-t_{m-1}=\delta''_{i+2}$ for all $m=1,\ldots,k$, where $\delta''_{i+2}$
is chosen so that
%
\begin{equation}\label{1213}
\qquad \bigl|L^{y_j}_{t_{m+1}}(\omega)-L^{y_j}_{t_{m-1}}(
\omega)\bigr|\le\frac1{2^{i+2}} \qquad\forall j=1,\ldots,n,\ \forall m=1,\ldots,k-1.
\end{equation}
Let $s_1,s_2\in[0,T]$, and assume that $s_1\le s_2$ and that
$s_2-s_1\le
\delta''_{i+2}$. There exists an $1\le m\le k-1$, such that
\[
t_{m-1}\le s_1\le s_2\le t_{m+1}.
\]
If $y,z\in D$ satisfy $\rho(y,z)\le\delta_{i+2}$, we
can find a $y_j\in Y$ such that $y\in B_{\rho}(y_j,\delta_{i+2})$.
If, in addition,
$L^{y}_{s_2}(\omega)\ge L^{z}_{s_1}(\omega)$, we have
%
\begin{eqnarray}\label{1214}
0&\leq& L^{y}_{s_2}(\omega)-L^{z}_{s_1}(
\omega)
\nonumber\\
& \le& L^{y}_{t_{m+1}}(\omega)-L^{z}_{t_{m-1}}(
\omega)
\nonumber
\\[-8pt]
\\[-8pt]
\nonumber
& \le&\bigl|L^{y}_{t_{m+1}}(\omega)-L^{y_j}_{t_{m+1}}(
\omega)\bigr| + \bigl|L^{y_j}_{t_{m+1}}(\omega)-L^{y_j}_{t_{m-1}}(
\omega)\bigr|
\\
&&{} +\bigl|L^{y_j}_{t_{m-1}}(\omega)-L^{y}_{t_{m-1}}(
\omega)\bigr| +\bigl|L^{y}_{t_{m-1}}(\omega)-L^{z}_{t_{m-1}}(
\omega)\bigr|,\nonumber
\end{eqnarray}
where the second inequality uses the fact that local time is
nondecreasing in $t$.
The second term to the right of the last inequality in (\ref{1214})
is less
than or equal to $2^{-(i+2)}$ by (\ref{1213}). The other three terms
are also
less than or equal to $2^{-(i+2)}$ by (\ref{123}) since
$\rho(y,y_j)\le\delta_{i+2}$
and $\rho(y,z)\le\delta_{i+2}$. Taking $\delta'_i=\delta''_{i+2}\land\delta_{i+2}$,
we get~(\ref{124}) on the larger set $[0,T']\times
D$ for some $T'\ge T$. Obviously this implies~(\ref{124}) as stated
in the case when $L^{y}_{s_2}(\omega)\ge L^{z}_{s_1}(\omega)$. A similar
argument gives~(\ref{124}) when $L^{y}_{s_2}(\omega)\le
L^{z}_{s_1}(\omega)$. Thus
(\ref{124}) is established.

 In what follows, we say that a function is locally uniformly
continuous
on a measurable set A in a locally compact metric space if it is uniformly
continuous on $A\cap K$
for all compact subsets $K\subseteq S$. Let $K_n$ be a sequence of compact
subsets of $S$ such that $S=\bigcup_{n=1}^\infty K_n$, and let $D'$ be a
countable dense subset of $S$. Let
\[
\hat\Om=\bigl\{\omega\mid L^y_t(\omega) \mbox{ is
locally uniformly continuous on } [0,\ze)\times D'\bigr\}.
\]
Let $ Q$ denote the rational numbers. Then
%
\begin{eqnarray}
\hat\Om^c&=&\mathop{\bigcup_{ s\in Q}}_{1\le n\le\infty} \bigl\{ \omega
\mid L^y_t(\omega) \mbox { is not uniformly continuous on}
\nonumber
\\[-8pt]
\\[-8pt]
\nonumber
&&\hspace*{85pt}\qquad [0,s]\times\bigl(K_n\cap D'\bigr);s<\ze\bigr\}.
\end{eqnarray}
Since $h_{x}>0$, it follows from (\ref{124}) and (\ref{sp2}) that
$P^x(\hat\Om^c)=0$ for all $x\in S$, or equivalently, that
%
\begin{equation}
P^x(\hat\Om)=1 \qquad \forall x\in S\label{1215} .
\end{equation}
We now construct a stochastic process
$\hat L=\{\hat L^y_t,(t,y)\in R_+\times S\}$
which is continuous on $[0,\ze)\times S$ and which is a version
of $L$. For $\omega\in\hat\Om$, let
$\{\tilde L^y_t(\omega),(t,y)\in[0,\ze)\times S\}$ be the continuous
extension of $\{L^y_t(\omega),(t,y)\in[0,\ze)\times D'\}$ to
$[0,\ze)\times S$. Set
%
\begin{eqnarray}
\hat L^y_t(\omega)&=&\tilde L^y_t(
\omega) \qquad\mbox{if }t<\ze(\omega)\label{1217a},
\\
 \hat L^y_t(\omega)&=&\mathop{\liminf_{s\uparrow\ze(\omega
)}}_{s\in Q} \tilde
L^y_t(\omega) \qquad\mbox{if }t\ge\ze(\omega) \label{1217}
\end{eqnarray}
and for $\omega\in\hat\Om^c$, set
\[
\hat L^y_t(\omega)\equiv0 \qquad\forall t,y\in R_+\times S.
\]
The stochastic
process $\{\hat L^y_t,(t,y)\in R_+\times S\}$ is well defined
and, clearly, is jointly continuous on $[0,\ze)\times S$.

We
now show that $\hat L$ is a local time by showing that for each $x,y\in S$,
%
\begin{equation}
\hat L^y_t=L^y_t\qquad \forall t\in
R_+, P^x \mbox{ almost surely.}\label{ltd}
\end{equation}
Recall that for each $z\in D'$, $\{L^z_t,
t\in R_+\}$ is increasing, $P^{x }$ almost surely. Hence, the same is
true for
$\{\tilde L^y_t,t<\ze\}$, and so the limit inferior in
(\ref{1217})
is actually a limit, $P^{x }$ almost surely.
Thus $\{\hat L^y_t,t\in R_+\}$ is
continuous and constant for $t\ge\ze$, $P^{x }$ almost surely.
Similarly, $L^y_t$, the local time for $X$ at $y$, is, by definition,
continuous in $t$ and constant for $t\ge\ze$, $P^{x }$ almost surely.
Now let us note that we could just as well have obtained
(\ref{124}) with $D'$ replaced by $D'\cup\{y\}$
and hence obtained (\ref{1215}) with $D'$ replaced by
$D'\cup\{y\}$ in the definition of $\hat\Om$.
Therefore if we take a
sequence $\{y_i\}_{i=1}^\infty$ with $y_i\in D'$ such that
$\lim_{i\to\infty}y_i=y$, we have that
%
\begin{equation}
\lim_{i\to\infty}L^{y_i}_t=L^y_t\qquad
\mbox{locally uniformly on $[0,\ze)$, $P^{x }$ a.s.}
\end{equation}
By the definition of $\hat L$, we also have
%
\begin{equation}
\lim_{i\to\infty}L^{y_i}_t=\hat L^y_t\qquad
\mbox{locally uniformly on $[0,\ze)$, $P^{x }$ a.s.}
\end{equation}
This shows that
%
\begin{equation}
\hat L^y_t= L^y_t\qquad  \forall t<
\ze, P^{x } \mbox{ a.s. }
\end{equation}
Since $\hat L^y_t$ and $ L^y_t$ are continuous in $t$ and constant
for $t\ge\ze$, we get (\ref{ltd}). This completes the proof of
Theorem~\ref{theo-jc} when $X$ is a transient Borel right process.

Now let $X$ be a recurrent Borel right process with state space $S$ and
continuous, strictly positive
$1$-potential densities
$u^{1} (x,y)$.\vadjust{\goodbreak}
Let $Y$ be the Borel right process obtained by killing $X$ at an independent
exponential time $\la$ with mean one. The 0-potential densities for
$Y$ are the 1-potential densities
for $X$. Thus we have a transient
Borel right process $Y$ with
continuous, strictly positive 0-potential densities $u^1(x,y)$. It is
easy to see that $L^y_{t\land\la}$ is a local time for $Y$.
Therefore, by what we have just shown for transient processes,
$L^y_t$ is continuous on
$S\times[0, \la)$, $P^{x}\times\nu$ almost surely,
where $\nu$ is the
probability measure of $\la$.
It now follows by Fubini's theorem that $L^y_t$ is continuous
$[0, q_i)\times S$ for
all $q_i\in Q$, $P^{x} $ almost surely, where~$Q$ is a countable dense
subset of $R_+$. This gives the proof when $X$ is recurrent.
\end{pf*}

We also can give a good local modulus of continuity for the local times.

\begin{theorem}\label{theo-62} Let $X =
(\Om, X_t, P^x
)$ be a Borel right process that satisfies all the hypotheses in
Theorem~\ref{theo-jc}. Let $d_{1}$ be a continuous metric or
pseudometric that dominates $d$ on $S\times S$. Fix $x_{0}\in S $, let
$T_{\delta}$ and $\mu_{\delta}$ be as in Lemma~\ref{lem-22q} and
assume that
(\ref{227}) holds. Then for almost every $t$,
%
\begin{equation}
\lim_{\delta\to0} \sup_{ d_{1}(x,x_{0})\le\delta/2 }\frac{
|L^{x}_{t}-L^{x_{0}}_{t}|}{\ov H_{T_{\delta},d_{1},\mu_{\delta}}
(\delta/4)}\le C
\bigl(L^{x_{0}}_{t}\bigr)^{1/2} \qquad\mbox{a.s.},
\label{21wwff}
\end{equation}
where $ \ov H_{T_{\delta},d_{1},\mu_{\delta}}(\delta/4)$ is
given in (\ref{220aa}).
\end{theorem}

\begin{pf} Let $\la$ be an independent mean one exponential. Note that for
a continuous function, the sup over any set can be evaluated by taking
the sup over a countable dense subset. Then using Theorem \ref
{theo-ke} and (\ref{21wwdd}),
\begin{eqnarray*}
&&\lim_{\delta\to0} \sup_{ d_{1}(x,x_{0})\le\delta/2
}\frac{ |L^{x}_{\la}-L^{x_{0}}_{\la}|}{\ov H_{T_{\delta},d_{1},\mu
_{\delta}} (\delta/4)}
\\
&&\qquad\leq\lim_{\delta\to0} \sup_{ d_{1}(x,x_{0})\le\delta/2
}\frac{ |L^{x}_{\la}+{ \th_{x} }/2-(L^{x_{0}}_{\la}+{
\th_{x_{ 0}} }/2)|}{\ov H_{T_{\delta},d_{1},\mu_{\delta}} (\delta/4)}\\
&&\qquad\quad{} +
\lim_{\delta\to0} \sup_{ d_{1}(x,x_{0})\le\delta/2 }\frac{ |
{ \th_{x} }/2- { \th_{x_{ 0}} }/2|}{\ov H_{T_{\delta
},d_{1},\mu_{\delta}} (\delta/4)}
\\
&&\qquad\leq C \biggl(L^{x_{0}}_{\la}+\frac{ \th_{x_{ 0}} }2
\biggr)^{1/2}+C \biggl(\frac{ \th_{x_{ 0}} }2 \biggr)^{1/2}\qquad
\mbox{a.s.},
\end{eqnarray*}
with respect to the product measure $ E^{x/h_{x}}E_{ \th} $.
Since $\th_{x_{ 0}} $ is the square of a normal random variable, for
any $\ep>0$ we have that $P_{\th} ( \th_{x_{ 0}} \leq\ep
)>0$. It then follows by Fubini's theorem that
%
\begin{equation}\quad
\lim_{\delta\to0} \sup_{ d_{1}(x,x_{0})\le\delta/2 }\frac{
|L^{x}_{\la}-L^{x_{0}}_{\la}|}{\ov H_{T_{\delta},d_{1},\mu
_{\delta}} (\delta/4)} \leq C
\bigl(L^{x_{0}}_{\la}+\ep \bigr)^{1/2}+C
\ep^{1/2} \qquad\mbox{a.s.} \label{lmodtt}
\end{equation}
The theorem follows by taking $\ep\rar0$ and then using Fubini's
theorem as in the last paragraph of the preceding proof.\vadjust{\goodbreak}
\end{pf}

\begin{pf*}{Proof of Theorem~\ref{lt-um}} We show below that for
any $\ep>0$, we can find $\gamma>0$
such that for all $x_{0}\in K$,
%
\begin{equation}
P \Bigl( \mathop{\sup_{x \in K}}_{ d(x_{0},x)\le\gamma}\th_{x
}^{1/2} \le\ep \Bigr)>0.
\label{pum1}
\end{equation}
The same proof leading to (\ref{lmodtt}), but using Theorem \ref
{theo-11a}, shows that
%
\begin{equation}
\lim_{\delta\to0} \mathop{\sup_{x,y\in K\cap
B_{d}(x_{0},\gamma)}}_{ d (x,y)\le\delta}\frac{|L^{x}_{\la}
-L^{y}_{\la}|}{J_{ d} ( d(x,y) /2) }\le C \Bigl(
\sup_{x\in
S}L^{x}_{\la}+\ep^{2}
\Bigr)^{1/2}\qquad \mbox{a.s.}\label{pum3}
\end{equation}
Using the compactness of $K$ this leads to
%
\begin{equation}
\lim_{\delta\to0} \mathop{\sup_{x,y\in K}}_{ d (x,y)\le\delta
}\frac{|L^{x}_{\la} -L^{y}_{\la}|}{J_{ d} ( d(x,y) /2) }\le C \Bigl(
\sup_{x\in S}L^{x}_{\la}+\ep^{2}
\Bigr)^{1/2}\qquad \mbox{a.s.}\label{pum4}
\end{equation}
The theorem follows by taking $\ep\rar0$ and then using Fubini's
theorem as in the previous proof.

Let $\Ga=\sup_{x\in K}u^{1}(x,x)$ and $\eta$ be a standard normal
random variable. For any $\ep>0$, we can find $\ep'>0$ such that
%
\begin{equation}
P \bigl(\Ga^{1/2} |\eta| \le\ep/2 \bigr)\ge2\ep'
\label{pum5}.
\end{equation}
Recalling Lemma~\ref{lem-1}, it follows that
%
\begin{equation}
\sup_{x\in K} P \bigl( \th_{x}^{1/2} \le\ep/2 \bigr)
\ge2\ep'. \label{636}
\end{equation}
By (\ref{21swb}), for some $\gamma'>0$, sufficiently small
%
\begin{equation}
P \biggl( \mathop{\sup_{ x,y\in K}}_{ d(x,y)\le\gamma'}\frac{
|\theta^{1/2}_{x}-\theta^{1/2}_{y}|}{J_{ d}(d(x ,y)/2)}\le30 \biggr)\ge1-
\ep'.\label{21swbqx}
\end{equation}
Under the hypothesis (\ref{18jv}), there exists a $0<\gamma\leq
\gamma'$, such that
%
\begin{equation}
P \biggl( \mathop{\sup_{ x,y\in K}}_{ d(x,y)\le\gamma} |\theta^{1/2}_{x }-
\theta^{1/2}_{y}| \le\frac{\ep}{2} \biggr)\ge1-
\ep'.\label{21swbq}
\end{equation}
For any $x_{0}\in K$, (\ref{pum1}) follows by taking
%
\begin{equation}
\th_{x}^{1/2}\le\th_{x_{0} }^{1/2} +\bigl|
\th_{x}^{1/2} - \th_{x_{0}}^{1/2} \bigr|
\end{equation}
and using (\ref{636}) and (\ref{21swbq}).
\end{pf*}

\section{\texorpdfstring{Further considerations of Theorem \protect\ref{theo-jc}}{Further considerations of Theorem 1.2}}
\label{sec-comp}

It is clear that Theorem~\ref{theo-jc} holds if $d$ in (\ref{18jj})
is replaced by a metric that dominates it. We use this observation to
show that Theorem~\ref{theo-jc}
gives the continuity results in~\cite{EK1}, Theorem 1.1.\vadjust{\goodbreak}

Let $X$ be a recurrent Borel right process with state space $S$ and
strictly positive $\al$-potential densities with respect to some
reference measure. Let $0$ be a distinguished point in $S$, and let
$u_{T_{0}}(x,y)$ denote the potential densities of the Borel right
process $Y$, which is $X$ killed the first time it hits $0$. In \cite
{EK1}, the authors show that when $X$ has a dual Borel right process,
$u_{T_{0}}(x,y)+u_{T_{0}}(y,x)$ is positive definite, so that
%
\begin{equation}
\ka(x,y)= \bigl(u_{T_{0}}(x,x)+u_{T_{0}}(y,y)-u_{T_{0}}(x,y)-u_{T_{0}}(y,x)
\bigr)^{1/2}\label{62}
\end{equation}
is a metric on $S$. In~\cite{EK1}, Theorem 1.1, they show that if for
every compact set $K\subseteq S$, there exists a probability measure
$\mu_{K}$ on $K$, such that
%
\begin{equation}
\lim_{\delta\to0} J_{K, \ka,\mu_{K}}( \delta) =0,\label{18qw}
\end{equation}
then the local times of $X$ are jointly continuous.

To see how this result follow from Theorem~\ref{theo-jc} let $\{
L_{t}^{y};(y,t)\in S\times R_{+}\}$ denote the local times of $X$.
Let $\tau(t)=\inf\{s\geq0 | L_{s}^{0}>t\}$ be the inverse local
time at $0$ and let $\la$ be an independent exponential random
variable with mean 1. Let $u_{\tau(\la)}(x,y)$ denote the potential
densities for the Borel right process $Z$, which is $X$ killed at $\tau
(\la)$.
It follows from~\cite{book}, (3.193), that
%
\begin{equation}
u_{\tau(\la)}(x,y)=u_{T_{0}}(x,y)+1.\label{63}
\end{equation}
Let $d(x,y)$ be the function defined in (\ref{12}) for the kernel
$u_{\tau(\la)}(x,y)$.

We now note that since $X$ has a dual Borel right process, so does $Y$.
Therefore $u_{T_{0}}(x,y)$, the potential of $Y$, satisfies (\ref
{49}). By (\ref{63}), $u_{\tau(\la)}(x,y)$ also satisfies (\ref
{49}) and, obviously, the $d_{2}$ metric for $u_{\tau(\la)}(x,y)$
[defined in (\ref{410})] is equal to $\ka(x,y)$.
Therefore, by (\ref{411}),
%
\begin{equation}
\tfrac{\sqrt{3}}{ 8} \,d (x,y)\le\ka(x,y) ,
\end{equation}
and consequently (\ref{18qw}) implies (\ref{18jj}).

Therefore, it follows from Theorem~\ref{theo-jc}, that $X$ has
continuous local times on $S\times[0,\tau(\la))$. Using Fubini's
theorem, as in the last paragraph of the proof of Theorem \ref
{theo-11a}, and the fact that $\lim_{t\rar\infty} \tau(t)=\infty
$, we see that
$X$ has jointly continuous local times on $S \times[0,\infty) $.

%


\printaddresses

\end{document}